\documentclass{article}

\usepackage{arxiv}

\usepackage[T1]{fontenc}    
\usepackage{hyperref}       
\usepackage{url}            
\usepackage{CJK}
\usepackage{booktabs}       
\usepackage{amsfonts}       
\usepackage{nicefrac}       
\usepackage{microtype}      
\usepackage{lipsum}
\usepackage{graphicx}
\usepackage{amsmath}
\usepackage[all]{xy}
\usepackage{amsthm}
\usepackage{eucal}
\usepackage{framed}

\newtheorem{theorem}{Theorem}
\theoremstyle{definition}
\newtheorem{definition}[theorem]{Definition}
\newtheorem{example}[theorem]{Example}

\newtheorem{lemma}{Lemma}
\newtheorem{remark}{Remark}

\newtheorem{proposition}{Proposition}
\newtheorem*{pf}{Proof}

\title{On Some Integral Representation Of $\zeta(n)$  Involving Nielsen's Generalized Polylogarithms And The Related Partition Problem}

\author{
  Xiaowei Wang(Potsdam)\thanks{This paper is written in November 2019}
   \\
}

\begin{document}
\maketitle

\begin{abstract}
In this paper, we study a family of single variable integral representations for some products of $\zeta(2n+1)$, where $\zeta(z)$ is Riemann zeta function and $n$ is positive integer. Such representation involves the integral $Lz(a,b):=\frac{1}{(a-1)!b!}\int_{0}^{1}\log^a (t)\log^b (1-t)dt/t$ with positive integers $a,b$, which is related to Nielsen's generalized polylogarithms. By analyzing the related partition problem, we discuss the structure of such integral representation, especially the condition of expressing products of $\zeta(2n+1)$ by finite $\mathbb{Q}(\pi)$-linear combination of $Lz(a,b)$.
\end{abstract}

\keywords{Riemann zeta function at integers \and integral representation \and Nielsen's generalized polylogarithms}
\section{Introduction}
It's well known that many number theoretic properties of $\zeta(2n+1)$ are nowadays still unsolved mysteries, such as the rationality (only known $\zeta(3)$ is irrational), transcendence and existence of closed-form functional equation that satisfied by $\zeta(2n+1)$. Thanks to the basic functional relation between gamma function $\Gamma(z)$ and sine function $\sin(z)$, i.e. $\Gamma(z)\Gamma(1-z)=\frac{\pi z}{\sin(\pi z)}=\pi z\csc(\pi z)$, one can explicitly express $\zeta(2n)$ by $r_{2n}\pi^{2n}$, where $r_{2n}$ is some rational number related to Bernoulli number $B_{2n}$. Unfortunately, to find such a simple analogous formula for $\zeta(2n+1)$ is considered to be impossible.
Studying the integral and series representations for $\zeta(2n+1)$ is somehow an important way to analyze the number theoretic properties of $\zeta(2n+1)$, for instance, F. Beukers' work \cite{1} is an excellent example. In this paper we discuss a class of integral representation with\\
\begin{equation*}
  Lz(a,b):=\frac{1}{(a-1)!b!}\int_{0}^{1}\frac{\log^{a-1}(t)\log^{b} (1-t)}{t}dt
\end{equation*}
Theoretically, many polynomials of $\zeta(n)$ can be represented by this integral. In fact, only $\zeta(2n+1)$ or $\zeta(2n+1)^d$ are interesting. Among all family of single variable integral representations that can represent polynomials of $\zeta(n)$, $Lz(a,b)$ is likely the most simple one. Via establishing some linear combination of $Lz(a,b)$ on $\mathbb{Q}(\pi)$, we are even able to express some $\zeta(2n+1)^d$. However, this method is in somehow restricted, which we shall discuss in the last section.\\
In fact, $(-1)^{a+b-1}Lz(a,b)$ is exactly a special value $S_{a,b}(1)$ of Nielsen's generalized polylogarithm $S_{a,b}(z)$, which was introduced by N. Nielsen\cite{2}.\\
\begin{equation*}
 S_{a,b}(z)=(-1)^{a+b-1}\frac{1}{(a-1)!b!}\int_{0}^{1}\frac{\log^a (t)\log^b (1-zt)}{t}dt
\end{equation*}
In most cases, this function is known for mathematical physicists in the context of quantum electrodynamics. Only few literatures \cite{3}\cite{4}\cite{5} concerned about the special case $S_{a,b}(1)$. However, number theoretic properties of $S_{a,b}(z)$ have seldom been studied.\\

Throughout the paper, integrals $\int_{0}^{1}f(t)dt$ may be regarded as improper, namely $\int_{0^{+}}^{1^{-}}f(t)dt$. $\zeta(s)$ denotes Riemann zeta function $\zeta(s):=\sum_{k=1}^{\infty}\frac{1}{k^{s}}$. In general, $N,m,n,k,i,j,l,a,b$ denote nonnegative integers. \\

\section{Preliminaries}
Our main result is based on the following well-known formula(\cite{6}, p45).
\begin{theorem}
For $z\in \mathbb{C}$ and $|z|<1$, we have\\
\begin{equation*}
  \Gamma(1+z)=\exp(-\gamma z+\sum_{k=2}^{\infty}(-1)^k \frac{\zeta(k)}{k}z^k)
\end{equation*}
\end{theorem}
\begin{pf}
By the definition of Gamma function\\
\begin{equation*}
  \Gamma(z)=\frac{1}{z}e^{-\gamma z}\prod_{n=1}^{\infty}(1+\frac{z}{n})^{-1}e^{z/n}
\end{equation*}
we can rewrite it as\\

\begin{equation*}
  \log\Gamma(1+z)=-\gamma z-\sum_{n=1}^{\infty} \log(1+\frac{z}{n})-\frac{z}{n}
\end{equation*}
When $|z|<1$, then $|z/n|<1$ for all $n=1,2,...$, thus we have\\
\begin{equation*}
\log(1+\frac{z}{n})=\sum_{k=1}^{\infty}(-1)^{k+1}\frac{z^k}{k n^k}=\frac{z}{n}+\sum_{k=2}^{\infty}(-1)^{k+1}\frac{z^{k}}{k n^k}
\end{equation*}
Therefore\\
\begin{align*}
  \log\Gamma(1+z)&=-\gamma z-\sum_{n=1}^{\infty} \sum_{k=2}^{\infty}(-1)^{k+1}\frac{z^{k}}{k n^k}+\frac{z}{n}-\frac{z}{n}\\
  &=-\gamma z+\sum_{n=1}^{\infty} \sum_{k=2}^{\infty}(-1)^{k}\frac{z^{k}}{k n^k}\\
  &=-\gamma z+\sum_{k=2}^{\infty}(-1)^{k}\frac{z^{k}}{k}\sum_{n=1}^{\infty}\frac{1}{n^k}\\
  &=-\gamma z+\sum_{k=2}^{\infty}(-1)^{k}\frac{\zeta(k)}{k}z^k
\end{align*}
Taking $\exp$ on both sides we get what we need to prove. The validity of changing the order of double sum is based on the normal convergence of $\Gamma(1+z)$.
\begin{flushright}
  $\Box$
\end{flushright}
\end{pf}
About the partition problem, the related notations we adopt are following.\\

For any positive integer $N$, a partition of $N$ is a way written $N$ into the sum of positive integers. Two sums that differ only in the order of their summands are considered the same partition. For given $N$, each partition of $N$ can be regarded as a finite multiset $\mathcal{M}=([n],\mu_{M})$ in which the underlying set is $[n]=\{n\in \mathbb{Z}:1\leq n\leq N\}$. Therefore $\mathcal{M}$ is determined by the multiplicity function $\mu_{M}:[n]\rightarrow \mathbb{Z}_{\geq 0}$ that satisfies\\
\begin{equation*}
\sum_{n=1}n\mu_{M}(n)=N
\end{equation*}
Now let $X=(x_1,...,x_{N})$, where $x_{n}=\mu_{M}(n)$, then $X$ totally determines $\mathcal{M}=([n],\mu_{M})$. Therefore the alternative way to define the partition of $N$ is by\\
\begin{definition}
\begin{equation*}
  \mathcal{P}(N):=\{X=(x_{1},x_{2},...,x_{N})\in \mathbb{Z}^{N}: \sum_{n=1}^{N}nx_{n}=N, 0\leq x_{n}\leq N \}
\end{equation*}
$\mathcal{P}(N)$ is called the partition set of $N$. Its element is called a partition element of $N$, which denoted by $X$.\\
\end{definition}
For given $X\in \mathcal{P}_{s}^{t}(N)$, the support and the norm of $X=(x_{1},x_2,...,x_{N})$ are defined by\\
\begin{definition}
\begin{align*}
 &Supp(X):=\{(n,x_n) :x_n>0\}\\
&\|X\|:=\sum_{n=1}^{N}x_{n}
\end{align*}
\end{definition}
The set of restricted partition of $N$ that has exactly $t$ parts and the size of each part is not less than $s(s>1)$, is denoted by $\mathcal{P}_{s}^{t}(N)$, namely \\
\begin{definition}
\begin{equation*}
  \mathcal{P}_{s}^{t}(N):=\{X=(0,...,0,x_{s},...,x_{N})\in \mathcal{P}(N):  \|X\|=t \}
\end{equation*}
\end{definition}
On could similarly define $\mathcal{P}_{s}(N), \mathcal{P}_{s}^{\geq t}(N), \mathcal{P}_{s}^{\leq t}(N)$ and so on. In this paper, $\mathcal{P}_{2}(N)$ is particularly more often used than others. That is\\
\begin{equation*}
  \mathcal{P}_{2}(N):=\{X=(0,x_{2},...,x_{N}): \sum_{n=2}^{N}nx_{n}=N ,\forall x_{n}\in\mathbb{Z} , 0\leq x_{n}\leq N\}
\end{equation*}

Further, the odd partition set and even partition set are defined by following
\begin{definition}
\begin{equation*}
 \mathcal{P}\mathcal{O}_{s}^{t}(N):=\{X=(x_{1},...,x_{N})\in\mathcal{P}_{s}^{t}(N): x_{2m}=0\text{ for all }m\in\mathbb{Z}\}
\end{equation*}
\begin{equation*}
 \mathcal{P}\mathcal{E}_{s}^{t}(N):=\{X=(x_{1},...,x_{N})\in\mathcal{P}_{s}^{t}(N):\forall x_{2m-1}=0\text{ for all }m\in\mathbb{Z} \}
\end{equation*}
\end{definition}
Note that for odd $N$, $\mathcal{P}\mathcal{E}_{s}^{t}(N)=\emptyset$ for all $t$, $\mathcal{P}\mathcal{O}_{s}^{t}(N)=\emptyset$ for all even $t$. Similarly, For even $N$, $\mathcal{P}\mathcal{O}_{s}^{t}(N)=\emptyset$ for all odd $t$.\\

Before discussing the relation between $Lz(a,b)$ and $\zeta(n)$, we shall introduce $Lz(a,b)$ and $lz(a,b)$.
\begin{definition}
For nonnegative integers $a,b$, define\\
\begin{equation*}
  lz(a,b):=\int_{0}^{1}log^a(t)log^{b}(1-t)dt
\end{equation*}
For positive integers $a,b$, define
\begin{equation*}
  Lz(a,b):=\frac{1}{a!b!}(lz(a,b)+alz(a-1,b)+blz(a,b-1))
\end{equation*}
\end{definition}
It's obvious to see that both $lz$ and $Lz$ are symmetric, namely $lz(a,b)=lz(b,a)$, $Lz(a,b)=Lz(b,a)$. In fact, we have\\

\begin{proposition}
\begin{equation*}
   Lz(a,b)=\frac{1}{(a-1)!b!}\int_{0}^{1}\frac{\log^{a-1}(t)\log^{b} (1-t)}{t}dt
\end{equation*}
or by the symmetry,
\begin{equation*}
   Lz(a,b)=\frac{1}{a!(b-1)!}\int_{0}^{1}\frac{\log^{b-1}(t)\log^{a} (1-t)}{t}dt
\end{equation*}
\end{proposition}
\begin{pf}
Only need to prove that\\
\begin{equation*}
  lz(a,b)+a lz(a-1,b)+b lz(a,b-1)=b\int_{0}^{1}\frac{\log^{b-1}(t)\log^{a} (1-t)}{t}dt
\end{equation*}
With integration by parts, one can see\\
\begin{align*}
  lz(a,b) & =-\int_{0}^{1}t(a\frac{1}{t}\log^{a-1}(t)\log^{b}(1-t)-b\frac{1}{1-t}\log^{a}(t)\log^{b-1}(1-t))dt \\
   & =-alz(a-1,b)+b\int_{0}^{1}\frac{t}{1-t}\log^{a}(t)\log^{b-1}(1-t)dt\\
\end{align*}
Still by substituting $x=1-t$, we obtain immediately\\
\begin{equation*}
  lz(a,b)=-alz(a-1,b)-blz(b-1,a)+b\int_{0}^{1}\frac{\log^{b-1}(x)\log^{a} (1-x)}{x}dx
\end{equation*}
That is what we need.
\begin{flushright}
  $\Box$
\end{flushright}
\end{pf}

\section{The relation between \texorpdfstring{$lz(a,b)$}{lz} and \texorpdfstring{$\zeta(n)$}{zeta(n)}}

Recall the simple relation between gamma function and beta function with $x,y\in \mathbb{C}$, $|x|,|y|<1$.\\
\begin{equation}\label{1}
f(x,y):=(1+x+y)B(1+x,1+y)=\frac{\Gamma(1+x)\Gamma(1+y)}{\Gamma(1+x+y)}
\end{equation}
Since $Re (1+x)>0, Re (1+y)>0$, the left hand side of the equation has the integral representation\\
\begin{equation*}
 B(1+x,1+y)=\int_{0}^{1}t^{x}(1-t)^{y}dt
\end{equation*}
Our strategy is follow: Applying Taylor's theorem for multivariate functions $f(x,y)$. On the one hand, at $(x,y)=(0,0)$ any all partial derivatives of $(1+x+y)\int_{0}^{1}t^{x}(1-t)^{y}dt$ can be evaluated explicitly. On the other hand, applying Theorem \ref{1} for $\Gamma(1+x)$, $\Gamma(1+y)$ and $\Gamma(1+x+y)$, then evaluate the expansion coefficients of $\frac{\Gamma(1+x)\Gamma(1+y)}{\Gamma(1+x+y)}$, namely\\
\begin{equation*}
  \frac{\Gamma(1+x)\Gamma(1+y)}{\Gamma(1+x+y)}=\exp(\sum_{n=2}^{\infty}(-1)^n \frac{\zeta(n)}{n}(x^n+y^n-(x+y)^n))
\end{equation*}
If we denote $(x+y)^n-x^n-y^n$ by $P_{n}(x,y)$ or $P_{n}$. Let\\
\begin{equation*}
D_{n,k}:=(-1)^{k(n+1)} \frac{\zeta(n)^{k}}{k!n^k}
\end{equation*}
Thus we can rewrite the above formula as\\
\begin{align*}
  \frac{\Gamma(1+x)\Gamma(1+y)}{\Gamma(1+x+y)}&=\prod_{n=2}^{\infty}\exp((-1)^{n+1}\frac{\zeta(n)}{n}P_{n})\\
  &=\prod_{n=2}^{\infty}(1+\sum_{k=1}^{\infty}(-1)^{k(n+1)}\frac{\zeta(n)^k}{n^k k!}P_{n}^k)\\
  &=\prod_{n=2}^{\infty}(1+\sum_{k=1}^{\infty}D_{n,k}P_{n}^k)\\
\end{align*}
Due to the normal convergence we can expand the last infinite product and rearrange terms with the order up to $\deg(P_{n}^{k})$, where $\deg(.)$ is the total degree of polynomial. Since $P_{n}(x,y)$ is homogeneous polynomial of $x,y$, therefore obviously $P_{n}^{k}(x,y)$ is also homogeneous polynomial of $x,y$ with $\deg(P_{n}^{k})=k\deg(P_{n})=nk$.\\

That is, we can expand and rearrange $\prod_{n=2}^{\infty}(1+\sum_{k=1}^{\infty}D_{n,k}P_{n}^k)$ as follow\\
\begin{align*}
  \prod_{n=2}^{\infty}(1+\sum_{k=1}^{\infty}D_{n,k}P_{n}^k)=&1+(D_{2,1}P_{2})+(D_{3,1}P_{3})\\
&+(D_{4,1}P_{4}+D_{2,2}P_{2}^{2})\\
&+(D_{5,1}P_{5}+D_{3,1}D_{2,1}P_{3}P_{2})\\
&+(D_{6,1}P_{6}+D_{4,1}D_{2,1}P_{4}P_{2}+D_{3,2}P_{3}^{2}+D_{2,3}P_{2}^{3})+...\\
\end{align*}
or\\
\begin{equation}\label{2}
  f(x,y)=\prod_{n=2}^{\infty}(1+\sum_{k=1}^{\infty}D_{n,k}P_{n}^k)=1+\sum_{N=2}^{\infty}\sum_{X\in \mathcal{P}_{2}(N)}\prod_{(\alpha,\beta)\in Supp(X)}D_{\alpha,\beta}P_{\alpha}^{\beta}
\end{equation}
Notice that\\
\begin{equation*}
  deg(\prod_{(\alpha,\beta)\in Supp(X)}D_{\alpha,\beta}P_{\alpha}^{\beta})=\sum_{(\alpha,\beta)\in Supp(X)}\alpha \beta=N
\end{equation*}
For fixed positive integers $a,b$, the term $x^a y^b$ has degree $a+b$. Therefore it only appears in $\sum_{X\in \mathcal{P}_{2}(a+b)}\prod_{(\alpha,\beta)\in Supp(X)}D_{\alpha,\beta}P_{\alpha}^{\beta}$. Now we can assume that\\
\begin{equation}\label{3}
 \sum_{X\in \mathcal{P}_{2}(N)}\prod_{(\alpha,\beta)\in Supp(X)}D_{\alpha,\beta}P_{\alpha}^{\beta}=\sum_{j=1}^{N-1}\rho_{N-j,j}x^{N-j}y^{j}
\end{equation}
Now we can evaluate $\rho_{a,b}$ in two ways. The first one is integral representation.\\
\begin{lemma}
For positive integers $a,b$, we have
\begin{equation*}
\rho_{a,b}=\frac{1}{a!b!}\frac{\partial^{a+b}}{\partial x^a \partial y^b}|_{(0,0)}f(x,y)=Lz(a,b)
\end{equation*}
\end{lemma}
\begin{pf}
Let $g(x,y)=1+x+y, h(x,y)=B(1+x,1+y)$, notice that $\frac{\partial g}{\partial x}=\frac{\partial g}{\partial y}=1$ for all $(x,y)$, therefore any one of second-order partial derivatives $\frac{\partial^2 g}{\partial x^2}, \frac{\partial^2 g}{\partial y^2}, \frac{\partial^2 g}{\partial x\partial y}$ vanishes. Using Leibniz rule for $x$-component\\
\begin{equation*}
\frac{\partial^{a}}{\partial x^a} gh=g \frac{\partial^{a}}{\partial x^a}h +a\frac{\partial^{a-1}}{\partial x^{a-1}}h
\end{equation*}
and then for $y$-component, we have
\begin{align*}
\frac{\partial^{a+b}}{\partial x^a \partial y^b}gh&=\frac{\partial^{b}}{\partial y^b}(g \frac{\partial^{a}}{\partial x^a}h +a\frac{\partial^{a-1}}{\partial x^{a-1}}h)\\
&=g\frac{\partial^{a+b}}{\partial x^a \partial y^b}h+b\frac{\partial^{a+b-1}}{\partial x^a \partial y^{b-1}}h+a\frac{\partial^{a+b-1}}{\partial x^{a-1} \partial y^{b}}h\\
\end{align*}
Its value at the point $(x,y)=(0,0)$ is\\
\begin{equation*}
\frac{\partial^{a+b}}{\partial x^a \partial y^b}|_{(0,0)}gh=(\frac{\partial^{a+b}}{\partial x^a \partial y^b}h+b\frac{\partial^{a+b-1}}{\partial x^a \partial y^{b-1}}h+a\frac{\partial^{a+b-1}}{\partial x^{a-1} \partial y^{b}}h)|_{(0,0)}
\end{equation*}
On the other hand notice that for positive integer $a,b$\\
\begin{align*}
\frac{\partial^{a+b}}{\partial x^a \partial y^b}|_{(0,0)}B(1+x,1+y)&=\int_{0}^{1}\frac{\partial^{a+b}}{\partial x^a \partial y^b}|_{(0,0)}t^x(1-t)^ydt\\
&=\int_{0}^{1}\log^a (t)\log^b (1-t)dt\\
&=lz(a,b)
\end{align*}
Therefore,\\
\begin{equation*}
\frac{\partial^{a+b}}{\partial x^a \partial y^b}|_{(0,0)}gh=lz(a,b)+a lz(a-1,b)+b lz(a,b-1)
\end{equation*}
Namely,\\
\begin{equation*}
\frac{1}{a!b!}\frac{\partial^{a+b}}{\partial x^a \partial y^b}|_{(0,0)}f(x,y)=Lz(a,b)
\end{equation*}
\begin{flushright}
  $\Box$
\end{flushright}
\end{pf}

Now we discuss the second approach.\\
\begin{lemma}
Assume that $n_{j}\geq 2, k_{j}\geq 1$, let
\begin{equation*}
  \prod_{j=1}^{K} P_{n_{j}}^{k_{j}}=\sum C_{\lambda,\mu}x^\lambda y^\mu
\end{equation*}
If $\lambda+\mu\neq\sum_{j=1}^{K}n_{j}k_{j}$ or $\lambda\mu=0$, then $C_{\lambda,\mu}=0$. Otherwise if $\mu=\sum_{j=1}^{K}n_{j}k_{j}-\lambda$ then $C_{\lambda,\mu}$ is given by\\
\begin{equation*}
C_{\lambda,\mu}=\sum_{\ell\in\mathfrak{S}(\mu)}\prod_{j=1}^{K}\prod_{i=1}^{k_{j}}\binom{n_{j}}{\ell_{ji}}
\end{equation*}
with $\mathfrak{S}(\mu)$ as follow\\
\begin{equation*}
  \mathfrak{S}(\mu)=\{\ell:\sum_{j=1}^{K}\sum_{i=1}^{k_{j}}\ell_{ji}=\mu, \forall\ell_{ji}\in \mathbb{Z}, 1\leq\ell_{ji}\leq n_{j}\}
\end{equation*}
\end{lemma}
\begin{pf}
Firstly, it is obviously that $P_{n_{j}}^{k_{j}}$ is homogeneous polynomial of degree $n_{j}k_{j}$ for each $j$ ,therefore $\prod_{j=1}^{K} P_{n_{j}}^{k_{j}}$ is also a homogeneous polynomial with $deg(\prod_{j=1}^{K} P_{n_{j}}^{k_{j}})=\sum_{j=1}^{K}n_{j}k_{j}$. On the other hand, notice that for all $n_{j}\geq2$, the coefficient of the terms $x^{n_{j}}$ and $y^{n_{j}}$ in $P_{n_{j}}$ are both $0$. Hence only $x^{\lambda}y^{\mu}$ with the conditions $\lambda+\mu=\sum_{j=1}^{K}n_{j}k_{j}$ and $\lambda\mu\neq0$ has nonzero coefficient.\\
Secondly, Assume that for each $j$ we have\\
\begin{equation*}
P_{n_{j}}^{k_j}=(\sum_{\ell_{ji}=1}^{n_j-1}\binom{n_{j}}{\ell_{ji}}x^{n_{j}-\ell_{ji}}y^{\ell_{ji}})^{k_{j}}
\end{equation*}
In this way, the coefficient of $x^{\lambda}y^{\mu}$ in the expansion of $\prod_{j=1}^{K} P_{n_{j}}^{k_{j}}$ should be the sum of all product of $\binom{n_{j}}{\ell_{ji}}$ that by choose $k_{j}$ coefficients from $P_{n_{j}}^{b_j}$ respectively and satisfying that\\
\begin{equation*}
  \sum_{j=1}^{K}\sum_{i=1}^{k_{j}}\ell_{ji}=\mu, \forall\ell_{ji}\in \mathbb{Z}, 1\leq\ell_{ji}\leq n_{j}
\end{equation*}
We denote such constraint by $\mathfrak{S}(\mu)$. In fact it is coincide with\\
\begin{equation*}
  \sum_{j=1}^{K}\sum_{i=1}^{k_{j}}n_{j}-\ell_{ji}=\lambda
\end{equation*}
since $\mu=\sum_{j=1}^{K}n_{j}k_{j}-\lambda$. Now for fixed $n_{j},k_{j},(1\leq j\leq K)$, $C$ only determined by $\lambda$ or $\mu$, we can form now on simplify this notation by $C_{\mu}$
\begin{flushright}
  $\Box$
\end{flushright}
\end{pf}
Above Lemma is for general integers $a_{j},b_{j}$, now we reformulate it and only aim to $Supp(X)$. Assume that $X\in \mathcal{P}_{2}(N)$, let\\
\begin{equation*}
  \prod_{(n,k)\in Supp(X)}P_{n}^{k}=\sum C_{b}(X)x^{a}y^{b}
\end{equation*}
Then if $b>N-\|X\|$ or $b<\|X\|$, then $\sum C_{b}(X)=0$. Otherwise, it is given by\\
\begin{equation*}
  \sum_{\ell\in\mathfrak{S}(b)}\prod_{(n.k)\in Supp(X)}\prod_{}^{k}\binom{n}{\ell}
\end{equation*}

Now we are able to represent $\zeta(N)$ and some $\prod_{\sum n_{j}=N}\zeta(n_{j})$ by $Lz(a,b)$ with $a+b=N$. Following we provide the representation structure of $Lz(a,b)$.\\
\begin{theorem}
(The Partition represented Relation) Assume that $a\leq b$, then
\begin{equation}\label{4}
  Lz(a,b)=\sum_{X\in \mathcal{P}_{2}(N),\|X\|\leq b}(c_{b}(X)\prod_{(n,k)\in supp(X)}\zeta(n)^{k})
\end{equation}
where $c_{X}$ is rational number related to $X$, and it can be evaluated by Lemma $2$
\end{theorem}
\begin{pf}
By the expansion (\ref{3}) and Lemma $1$, we have\\

\begin{align*}
  \sum_{j=1}^{N-1}Lz(N-j,j)x^{N-j}y^{j}&=\sum_{j=1}^{N-1}\rho_{N-j,j}x^{N-j}y^{j}\\
    &=\sum_{X\in \mathcal{P}_{2}(N)}\prod_{(n,k)\in supp(X)}D_{n,k}P_{n}^{k}\\
  &=\sum_{X\in \mathcal{P}_{2}(N)}\prod_{(n,k)\in supp(X)}D_{n,k}\prod_{(n,k)\in supp(X)}P_{n}^{k}\\
\end{align*}
It remains to show that the coefficient of $x^{N-b}y^{b}$ on the right hand side has the form (\ref{4})\\
On the other hand, by Lemma $2$, we notice that for any term $Q$ of $\prod_{(n,k)\in supp(X)}P_{n}^{k}$\\
\begin{equation*}
\deg_{y}(Q)\geq\sum_{(n,k)\in supp(X)}k=\|X\|
\end{equation*}
and\\
\begin{equation*}
\deg_{y}(Q)\leq \sum_{(n,k)\in supp(X)}(n-1)k=N-\|X\|
\end{equation*}
Therefore
\begin{equation*}
\prod_{(n,k)\in supp(X)}P_{n}^{k}=\sum_{j=\|X\|}^{N-\|X\|} C_{j}(X)x^{N-j}y^{j}
\end{equation*}
That is to say,
\begin{equation*}
\sum_{j=1}^{N-1}Lz(N-j,j)x^{N-j}y^{j}=\sum_{X\in \mathcal{P}_{2}(N)}q_{X}\sum_{j=\|X\|}^{N-\|X\|} C_{j}(X)x^{N-j}y^{j}
\end{equation*}
where $q_{X}=\prod_{(n,k)\in supp(X)}D_{n,k}$. Now compare the coefficients on both sides, we have\\
\begin{equation*}
Lz(N-b,b)=\sum_{X\in \mathcal{P}_{2}(N),\|X\|\leq b}C_{b}(X)q_{X}
\end{equation*}
Recalling that $D_{n,k}:=(-1)^{k(n+1)} \frac{\zeta(n)^{k}}{k!n^k}$, so there is rational number $\widetilde{C}(X)$ such that\\
\begin{equation*}
  q_{X}=\prod_{(n,k)\in supp(X)}D_{n,k}=\widetilde{C}(X)\prod_{(n,k)\in supp(X)}\zeta(n)^{k}
\end{equation*}
Therefore\\
\begin{equation*}
Lz(N-b,b)=\sum_{X\in \mathcal{P}_{2}(N),\|X\|\leq b} C_{b}(X)\widetilde{C}(X)\prod_{(n,k)\in supp(X)}\zeta(n)^{k}
\end{equation*}
Let $ C_{b}(X)\widetilde{C}(X)$ denoted by $c_{b}(X)$. It's obviously rational. We finally have\\
\begin{equation*}
  Lz(a,b)=\sum_{X\in \mathcal{P}_{2}(N),\|X\|\leq b}(c_{b}(X)\prod_{(n,k)\in supp(X)}\zeta(n)^{k})
\end{equation*}
\begin{flushright}
  $\Box$
\end{flushright}
\end{pf}

The expression of $c_{b}(X)$ is given by following\\
\begin{theorem}
Assume that $Supp(X)=\{(n_{j},k_{j}):1\leq j\leq J\}$, then
\begin{equation*}
  c_{b}(X)=(-1)^{N+\|X\|}\prod_{j=1}^{J}\frac{1}{k_{j}!n_{j}^{k_{j}}}\sum_{\ell\in \mathfrak{S}(b)}\prod_{i=1}^{k_{j}}\binom{n_{j}}{\ell_{ji}}
\end{equation*}
or rewrite as
\begin{equation*}
  c_{b}(X)=(-1)^{N+\|X\|}\prod_{(n,k)\in supp(X)} \frac{1}{k!}\sum_{\ell\in \mathfrak{S}(b)}\prod^{k}\frac{(n-1)!}{\ell!(n-\ell)!}
\end{equation*}
where $\mathfrak{S}(b)$ is given by\\
\begin{equation*}
  \mathfrak{S}(b)=\{\ell\in\mathbb{Z}: \sum_{j=1}^{K}\sum_{i=1}^{k_{j}}\ell_{ji},1\leq \ell\leq b-1 \}
\end{equation*}
\end{theorem}
\begin{pf}
The proof is straightforward. Recall that $ c_{b}(X)=C_{b}(X)\widetilde{C}(X)$, now on the one hand we have,\\
\begin{align*}
 \widetilde{C}(X)&=\prod_{(n,k)\in supp(X)}(-1)^{k(n+1)} \frac{1}{k!n^k}\\
 &=(-1)^{\sum _{j=1}^{J}k_{j}(n_{j}+1)}\prod_{j=1}^{J}\frac{1}{k_{j}!n_{j}^{k_{j}}}\\
 &=(-1)^{N+\|X\|}\prod_{j=1}^{J}\frac{1}{k_{j}!n_{j}^{k_{j}}}
\end{align*}
On the other hand, by Lemma $2$\\
\begin{equation*}
C_{b}(X)=\sum_{\sum \ell=b}\prod_{j=1}^{J}\prod_{i=1}^{k_{j}}\binom{n_{j}}{\ell_{ji}}
\end{equation*}
Multiply $C_{b}(X)$ and $\widetilde{C}(X)$ together, then we have what we need.
\begin{flushright}
  $\Box$
\end{flushright}
\end{pf}

\section{Properties of \texorpdfstring{$Lz(a,b)$}{Lz} and \texorpdfstring{$lz(a,b)$}{lz}}

\begin{proposition}
\begin{equation*}
 Lz(a,1)=\frac{(-1)^a}{a!}\int_{0}^{\infty}\frac{z^{a}}{e^z-1}dz=(-1)^a\zeta(a+1).
\end{equation*}
\end{proposition}
The proof is also easy, consider the substitution of $t=1-e^{-z}$. This formula connects $Lz(a,1)$ to the well-known formula about $\zeta(n)$.\\

\begin{proposition}
\begin{equation*}
 Lz(2n-1,2)=n\zeta(2n+1)-\sum_{j=2}^{n}\zeta(j)\zeta(2n-j+1)
\end{equation*}
\begin{equation*}
 Lz(2n-2,2)=(n-\frac{1}{2})\zeta(2n)-\sum_{j=2}^{n}\zeta(j)\zeta(2n-j)
\end{equation*}
or mixing them, as\\
\begin{equation*}
  Lz(a,2)=\frac{a+1}{2}\zeta(a+2)-\sum_{j=2}^{\lfloor\frac{a}{2}\rfloor+1}\zeta(j)\zeta(a+2-j)
\end{equation*}
\end{proposition}
\begin{pf}

By Theorem $7$, only need to consider the restricted partition that has merely one or two parts. For the case I. $N=2n+1$\\
For exactly one-part partition, there is only one element $X_{1}\in \mathcal{P}_{2}(N)$, and $Supp(X_{1})=\{(N,1)\}$. Similarly, for two-part partition, there are $n-1$ elements: $X_{2,j}\in \mathcal{P}_{2}(N)$ with $Supp(X_{2,j})=\{(j,1),(N-j,1)\}$, where $j=2,...,n$. Therefore\\
\begin{equation*}
  Lz(2n-1,2)=c_{2}(X_{1})\zeta(2n+1)+\sum_{j=2}^{n}c_{2}(X_{2,j})\zeta(j)\zeta(2n-j+1)
\end{equation*}
By Theorem $8$, it easy to get the coefficients \\
\begin{align*}
   & c_{2}(X_{1})=\binom{N}{2}\frac{1}{N}=n \\
   & c_{2}(X_{2,j})=-\binom{j}{1}\binom{N-j}{1}\frac{1}{j}\frac{1}{N-j}=-1
\end{align*}
Therefore\\
\begin{equation*}
 Lz(2n-1,2)=n\zeta(2n+1)-\sum_{j=2}^{n}\zeta(j)\zeta(2n-j+1)
\end{equation*}
case II. $N=2n$.\\
For exactly one-part partition, there is only one element $X_{1}\in \mathcal{P}_{2}(N)$, and $Supp(X_{1})=\{(N,1)\}$. Similarly, for two-part partition, there are $n-1$ elements: $X_{2,j}\in \mathcal{P}_{2}(N)$ with $Supp(X_{2,j})=\{(j,1),(N-j,1)\}$, where $j=2,...,n$. Therefore\\
\begin{equation*}
  Lz(2n-2,2)=c_{2}(X_{1})\zeta(2n)+\sum_{j=2}^{n}c_{2}(X_{2,j})\zeta(j)\zeta(2n-j)
\end{equation*}
By Theorem $8$, it easy to get the coefficients \\
\begin{align*}
   & c_{2}(X_{1})=\binom{N}{2}\frac{1}{N}=\frac{2n-1}{2} \\
   & c_{2}(X_{2,j})=-\binom{j}{1}\binom{N-j}{1}\frac{1}{j}\frac{1}{N-j}=-1
\end{align*}
Therefore\\
\begin{equation*}
 Lz(2n-2,2)=(n-\frac{1}{2})\zeta(2n)-\sum_{j=2}^{n}\zeta(j)\zeta(2n-j)
\end{equation*}
Let $a=2n-1$ or $a=2n-2$, rewrite those two formulas we have\\
\begin{equation*}
  Lz(a,2)=\frac{a+1}{2}\zeta(a+2)-\sum_{j=2}^{\lfloor\frac{a}{2}\rfloor+1}\zeta(j)\zeta(a+2-j)
\end{equation*}
\end{pf}

there is a relation between $Lz(a,b)$ and multiple zeta function. Following example shows that the symmetric property of $Lz(a,b)$ implies a nontrivial relation between $Lz(a,b)$ and multiple zeta function for argument of integers $\zeta(n_{1},n_{2})$. For larger $a,b$, we need more techniques. \\
\begin{example}
Consider $Lz(2,1)=Lz(1,2)$ . One the one hand, \\
\begin{align*}
  Lz(2,1)&=\int_{0}^{1}\frac{\log(t)\log(1-t)}{t}dt\\
  &=-\int_{0}^{1}\log(t)(1+\sum_{n=1}^{\infty}\frac{t^n}{n+1})dt\\
  &=-(\int_{0}^{1}\log(t)dt+\sum_{n=1}^{\infty}\frac{1}{n+1}\int_{0}^{1}t^n\log(t)dt)\\
  &=\sum_{n=0}^{\infty}\frac{1}{(n+1)^3}=\zeta(3)
\end{align*}
However, on the other hand, with the similar trick, for $|t|<1$, we have the expansion\\
\begin{equation*}
\frac{\log^2(1-t)}{t}=\sum_{n=1}^{\infty}S_{n+1}^{(2)}t^n
\end{equation*}
where $S_{n}^{(2)}$ denotes $\sum_{i+j=n,i,j\geq1}\frac{1}{ij}$. then\\\\
\begin{equation*}
Lz(1,2)=\frac{1}{2!}\int_{0}^{1}\frac{\log^2(1-t)}{t}dt=\frac{1}{2}\sum_{n=2}^{\infty}\frac{ S_{n}^{(2)}}{n}
\end{equation*}
Notice that for $S_{n}^{(2)}$\\
\begin{equation*}
  S_{n}^{(2)}=\frac{1}{n}\sum_{k=1}^{n-1}\frac{n}{k(n-k)}=\frac{1}{n}\sum_{k=1}^{n-1}(\frac{1}{k}+\frac{1}{n-k})=\frac{2}{n}\sum_{k=1}^{n-1}\frac{1}{k}
\end{equation*}
Therefore, above $Lz(1,2)$ can be reformulated as\\
\begin{equation*}
 Lz(1,2) =\frac{1}{2}\sum_{n=2}^{\infty}\frac{ S_{n}^{(2)}}{n}=\sum_{n>m\geq1}\frac{1}{n^2 m}=\zeta(2,1)
\end{equation*}
Hence we prove Euler identity $\zeta(3)=\zeta(2,1)$ by using $Lz(a,b)=Lz(b,a)$.
\end{example}
\begin{remark}
There is another proof using a series involving Striling numbers, see\cite{7}.
\end{remark}
In fact, generally we have\\
\begin{theorem}
(Symmetry of Series Representation) For any positive integer $a,b$, we have\\
\begin{equation*}
\frac{1}{a!}\sum_{n=a}^{\infty}\frac{S_{n}^{(a)}}{n^{b}} =\frac{1}{b!}\sum_{n=b}^{\infty}\frac{S_{n}^{(b)}}{n^{a}}
\end{equation*}
where\\
\begin{equation*}
  S_{n}^{(k)}=\sum_{\sum m_{j}=n}\prod_{j=1}^{k}m_{j}^{-1}, \text{     }  m_{j}\in\mathbb{Z}^{+}, j=1,...,k,
\end{equation*}
In particular, $ S_{n}^{(1)}=\frac{1}{n}$
\end{theorem}
\begin{pf}
Given positive integer $a,b$,\\
\begin{align*}
  Lz(a,b) & =\frac{1}{(a-1)!b!}\int_{0}^{1}\frac{\log^{a-1}(t)\log^{b}(1-t)}{t}dt \\
   &= \frac{1}{(a-1)!b!}\int_{0}^{1}\log^{a-1}(t)(-\sum_{n=1}^{\infty}\frac{t^n}{n})^b \frac{1}{t}dt\\
   &=\frac{(-1)^b}{(a-1)!b!}\int_{0}^{1}\log^{a-1}(t)\sum_{n=b}^{\infty}S_{n}^{(b)}t^{n-1}dt\\
   &=\frac{(-1)^b}{(a-1)!b!}\sum_{n=b}^{\infty}S_{n}^{(b)}\int_{0}^{1}\log^{a-1}(t)t^{n-1}dt\\
\end{align*}
By the substitution $t=e^{-z}$, it turns out to be\\
\begin{align*}
   Lz(a,b) &= \frac{(-1)^b}{(a-1)!b!}\sum_{n=b}^{\infty}S_{n}^{(b)}\int_{0}^{+\infty}(-z)^{a-1}e^{-nz}dz\\
   & =\frac{(-1)^{a+b-1}}{b!}\sum_{n=b}^{\infty}\frac{S_{n}^{(b)}}{n^{a}}
\end{align*}
On the other hand, by the similar method, we have\\
\begin{equation*}
  Lz(b,a)=\frac{(-1)^{a+b-1}}{a!}\sum_{n=a}^{\infty}\frac{S_{n}^{(a)}}{n^{b}}
\end{equation*}
Since $Lz(a,b)=Lz(b,a)$ holds for all positive integers $a,b$, therefore\\
\begin{equation*}
  \frac{(-1)^{a+b-1}}{b!}\sum_{n=b}^{\infty}\frac{S_{n}^{(b)}}{n^{a}}=\frac{(-1)^{a+b-1}}{a!}\sum_{n=a}^{\infty}\frac{S_{n}^{(a)}}{n^{b}}
\end{equation*}
namely\\
\begin{equation*}
 \frac{1}{a!}\sum_{n=a}^{\infty}\frac{S_{n}^{(a)}}{n^{b}} =\frac{1}{b!}\sum_{n=b}^{\infty}\frac{S_{n}^{(b)}}{n^{a}}
\end{equation*}
\end{pf}
A straightforward corollary is that, if $b=1$, we have\\
\begin{equation*}
  \zeta(a+1)=\frac{1}{a!}\sum_{n=a}^{\infty}\frac{S_{n}^{(a)}}{n}
\end{equation*}

On the other hand, by the substitution of $t=\sin^{2}(\theta)$ we have\\
\begin{proposition}
\begin{equation*}
  lz(a,b)=2^{a+b+1}\int_{0}^{\frac{\pi}{2}}\sin(\theta)\cos(\theta)\log^{a}\sin(\theta)\log^{b}\cos(\theta)d\theta
\end{equation*}
\begin{equation*}
  Lz(a,b)=\frac{2^{a+b}}{a!(b-1)!}\int_{0}^{\frac{\pi}{2}}\cot(\theta)\log^{b-1}\sin(\theta)\log^{a}\cos(\theta)d\theta
\end{equation*}
\end{proposition}

\section{Examples of Establishing the Integral Representations}
\begin{itemize}
  \item $N=3$
\end{itemize}
As the first example, for $N=3$, the integral representation is trivial. But as it already demonstrated in the last section, those two equivalent representations imply some nontrivial relation $\zeta(3)=\zeta(2,1)$. \\
\begin{equation*}
  \zeta(3)=Lz(1,2)=\frac{1}{2}\int_{0}^{1}\frac{\log^2(1-t)}{t}dt
\end{equation*}
\begin{equation*}
  \zeta(3)=Lz(2,1)=\int_{0}^{1}\frac{\log(t)\log(1-t)}{t}dt
\end{equation*}
\begin{itemize}
  \item $N=4$
\end{itemize}
\begin{align*}
  &\zeta(4)=-Lz(3,1)=-\frac{1}{2}\int_{0}^{1}\frac{\log^2(t)\log(1-t)}{t}dt\\
  &\zeta(4)=-4Lz(2,2)=-2\int_{0}^{1}\frac{\log(t)\log^2(1-t)}{t}dt\\
   &\zeta(4)=-Lz(1,3)=-\frac{1}{6}\int_{0}^{1}\frac{\log^3(1-t)}{t}dt
\end{align*}
They correspond to following series representations respectively. Let $S_{n}^{(3)}$ denotes $\sum_{n_{1}+n_{2}+n_{3}}\frac{1}{n_{1}n_{2}n_{3}}$, then \\
\begin{align*}
  &\zeta(4)=\sum_{n=1}^{\infty}\frac{1}{n^4}\\
  &\zeta(4)=2\sum_{n=2}^{\infty}\frac{S_{n}^{(2)}}{n^2}=4\zeta(3,1)\\
  &\zeta(4)=\frac{1}{6}\sum_{n=3}^{\infty}\frac{S_{n}^{(3)}}{n}
\end{align*}
Following we only concern about the integral representations.\\
\begin{itemize}
  \item $N=5$
\end{itemize}
\begin{align*}
  Lz(4,1) & =\zeta(5) \\
  Lz(3,2) & = 2\zeta(5)-\zeta(2)\zeta(3)
\end{align*}
By mixing them, we obtain\\
\begin{equation*}
  \zeta(3)=\frac{1}{\zeta(2)}(2Lz(4,1)-Lz(3,2))=\frac{1}{2\pi^2}\int_{0}^{1}\frac{\log^{2}(t)\log(1-t)}{t}\log\frac{t^4}{(1-t)^3}dt
\end{equation*}
This is a new nontrivial integral representation of Ap\'{e}ry's constant. On the other hand, we have another integral representation for $\zeta(5)$.\\
\begin{equation*}
  \zeta(5)=\frac{1}{8}\int_{0}^{1}\frac{\log^2(1-t)}{t}(\log^2 (t)+\frac{\pi^2}{3})dt
\end{equation*}
Once notice that\\
\begin{equation*}
  \frac{1}{8}\int_{0}^{1}\frac{\log^2(1-t)}{t}(\log (t)\frac{2\pi}{\sqrt{3}})dt=-\frac{\pi\zeta(4)}{8\sqrt{3}}
\end{equation*}
By plus to above formula, we obtain\\
\begin{equation*}
  \zeta(5)-\frac{\pi\zeta(4)}{8\sqrt{3}}=\frac{1}{8}\int_{0}^{1}\frac{\log^2(1-t)\log^2(e^{\pi/\sqrt{3}}t)}{t}dt
\end{equation*}
\begin{itemize}
  \item $N=6$
\end{itemize}
$\mathcal{P}_{2}(6)$ has $4$ elements:\\
\begin{equation*}
  X_{1}=(0,0,0,0,0,1), X_{2}=(0,1,0,1,0,0), X_{3}=(0,0,2,0,0,0), X_{4}=(0,3,0,0,0,0)
\end{equation*}
\begin{align*}
  &Supp(X_{1})=\{(6,1)\}\\
  &Supp(X_{2})=\{(2,1),(4,1)\}\\
  &Supp(X_{3})=\{(3,2)\}\\
  &Supp(X_{4})=\{(2,3)\}
\end{align*}
By (\ref{3})\\
\begin{equation*}
  Lz(5,1)x^5 y+Lz(4,2)x^4 y^2+Lz(3,3)x^3 y^3=D_{6,1}P_{6}+D_{2,1}D_{4,1}P_{2}P_{4}+D_{3,2}P_{3}^{2}+D_{2,3}P_{2}^3
\end{equation*}
Comparing the coefficients, we have\\
\begin{align*}
 & Lz(5,1)=-\zeta(6)\\
  &Lz(4,2) =\frac{1}{2}\zeta(3)^2+\zeta(2)\zeta(4)-\frac{5}{2}\zeta(6)=\frac{1}{2}\zeta(3)^2-\frac{\pi^6}{1260}\\
  &Lz(3,3)=\zeta(3)^2+\frac{3}{2}\zeta(2)\zeta(4)-\frac{10}{3}\zeta(6)-\frac{1}{6}\zeta(2)^3=\zeta(3)^2-\frac{23\pi^6}{15120}\\
\end{align*}
By mixing above two equations, we can rewrite a more interesting but sophisticated formula.\\
\begin{equation*}
  \int_{0}^{1}\frac{\log^{2}(t)\log^{2}(1-t)}{t}\log\frac{(1-t)^{12}}{t^{23}}dt=6\zeta(3)^2
\end{equation*}

\begin{itemize}
  \item $N=7$
\end{itemize}
\begin{align*}
   & Lz(6,1)=\zeta(7) \\
   & Lz(5,2)=3\zeta(7)-\zeta(2)\zeta(5)-\zeta(4)\zeta(3)\\
   & Lz(4,3)=5\zeta(7)-2\zeta(2)\zeta(5)-\frac{5}{4}\zeta(4)\zeta(3)
\end{align*}
reformulate them, we have interesting similar representation of $\zeta(3)$ and $\zeta(5)$.\\
\begin{align*}
   & \frac{3}{5}\zeta(2)\zeta(5)=Lz(6,1)+Lz(5,2)-\frac{4}{5}Lz(4,3) \\
   & \frac{3}{4}\zeta(4)\zeta(3)=Lz(6,1)-2Lz(5,2)+Lz(4,3)
\end{align*}

\begin{itemize}
  \item $N=8$
\end{itemize}
\begin{align*}
  &Lz(7,1)=-\zeta(8)\\
  &Lz(6,2)=\zeta(3)\zeta(5)-\frac{\pi^8}{7560}\\
  &Lz(5,3)=3\zeta(3)\zeta(5)-\frac{\pi^2 \zeta(3)^2}{12}-\frac{61\pi^8}{226800}\\
  &Lz(4,4)=4\zeta(3)\zeta(5)-\frac{\pi^2 \zeta(3)^2}{6}-\frac{499\pi^8}{1814400}
\end{align*}
\begin{itemize}
  \item $N=9$
\end{itemize}
\begin{align*}
  &Lz(8,1)=\zeta(9)\\
&Lz(7,2)=4\zeta(9)-\zeta(2)\zeta(7)-\zeta(4)\zeta(5)-\zeta(6)\zeta(3)\\
\end{align*}
\begin{align*}
  Lz(6,3)&=\frac{28}{3}\zeta(9)-3\zeta(2)\zeta(7)-\frac{7}{2}\zeta(4)\zeta(5)-\frac{7}{2}\zeta(6)\zeta(3)+\zeta(2)\zeta(3)\zeta(4)+\frac{1}{2}\zeta(2)^2\zeta(5)+\frac{1}{6}\zeta(3)^3\\
&=\frac{\zeta(3)^3}{6}+\frac{28\zeta(9)}{3}-\frac{\pi^6\zeta(3)}{540}-\frac{\pi^4\zeta(5)}{40}-\frac{\pi^2\zeta(7)}{2}\\
\end{align*}
\begin{align*}
  Lz(5,4)&=14\zeta(9)-5\zeta(2)\zeta(7)-6\zeta(4)\zeta(5)-\frac{35}{6}\zeta(6)\zeta(3)+\frac{5}{2}\zeta(2)\zeta(3)\zeta(4)+\zeta(2)^2\zeta(5)+\frac{1}{2}\zeta(3)^3-\frac{1}{6}\zeta(2)^3\zeta(3)\\
&=\frac{\zeta(3)^3}{2}+14\zeta(9)-\frac{\pi^6\zeta(3)}{432}-\frac{7\pi^4\zeta(5)}{180}-\frac{5\pi^2\zeta(7)}{6}\\
\end{align*}
When $N$ become larger, $a,b$ become closer, the expression of $Lz(a,b)$ would be more complicated. In fact, one can prove following statement\\
\begin{theorem}
If integer $N>20$, then there always exist partition $X\in \mathcal{P}_{2}(N)$, such that $\prod_{(n,k)\in Supp(X)}\zeta(n)^k$ cannot be represented by finite $\mathbb{Q}(\pi)$-linear combination of $Lz(a,b)$ for all $a,b\in \mathbb{Z}^+$ with $a+b\leq N$ by using the partition represented relation (Theorem $7$).\\
\end{theorem}
\begin{remark}
It's still unknown whether there is any functional equation in closed form that satisfied by $\zeta(n)$ for any different $n$, therefore such $\mathbb{Q}(\pi)$-linear combination of $Lz(a,b)$ may be constructed in other ways that differ from the partition represented relation. Hence all the \emph{representation}s in the following proof are referred to the representations by only using the partition represented relation.
\end{remark}
\begin{pf}
Firstly, for fixed $N$, let
\begin{equation*}
  \Pi:\mathcal{P}(N)\rightarrow \mathbb{R}; X\mapsto \prod_{(n,k)\in supp(X)}\zeta(n)^k
\end{equation*}
Assume that $X=X_{1}+X_{2}\in \mathcal{P}_{2}(N)$ with $X_{1}\in \mathcal{P}\mathcal{E}_{2}(N)$, $X_{2}\in \mathcal{P}\mathcal{O}_{3}(N)$. Since for all positive even number $2n$, $\zeta(2n)$ can be represented as $q_{n}\pi^{2n}$ with $q_{n}\in \mathbb{Q}$, then in other words $\Pi(X_{1})\in \mathbb{Q}(\pi)$. Conversely, if $\Pi(X_{2})$ cannot be represented as finite $\mathbb{Q}(\pi)$-linear combination of $Lz(a,b)$ for $a+b\leq N$, then $\Pi(X)$ neither.\\
Therefore, it remains to prove that there always exist $X\in \mathcal{P}\mathcal{O}_{3}(N)$ for sufficiently large $N$, such that $\Pi(X)$ cannot be represented by finite $\mathbb{Q}(\pi)$-linear combination of $Lz(a,b)$. For even or odd $N$, such $X$ is constructed by different way.
\begin{itemize}
  \item Case I. Suppose that $N=2M+1$. \\
\end{itemize}
Let\\
\begin{equation*}
  T(N)=\bigsqcup_{t=1}^{\infty}\mathcal{P}\mathcal{O}_{3}^{2t+1}(N)
\end{equation*}
Notice that for $Y\in \mathcal{P}_{2}(N)\backslash T(N)$, we can always assume that $\Pi(\widetilde{X})$ can be represented by finite $\mathbb{Q}(\pi)$-linear combination of $Lz(a,b)$ for $a+b< N$. Because otherwise, we come to a smaller $N_{1}<N$, and therefore we could repeat the process by starting from $N_{1}$ instead of $N$. That is. it reasonable to assume that if $Y\in \mathcal{P}_{2}(N)\backslash T(N)$, then
\begin{equation*}
  \sum Y=\sum p_{j} Lz(a_{j},b_{j})
\end{equation*}
with $a_{j}+b_{j}< N$. Now Suppose that \\
\begin{equation*}
\mathcal{P}\mathcal{O}_{3}^{2t+1}(N)=\{X_{1}^{(2t+1)},X_{2}^{(2t+1)},...,X_{r_{2t+1}}^{(2t+1)}\}
\end{equation*}
with $|\mathcal{P}\mathcal{O}_{3}^{2t+1}(N)|=r_{2t+1}$. According to theorem $7$, $\Pi(X_{i}^{(2t+1)})$ only appears in the representation formulas $Lz(N-2t-1,2t+1), Lz(N-2t-2,2t+2),...,Lz(2t+1,N-2t-1)$. In fact, due to the symmetric property of $Lz(a,b)$, only $Lz(N-2t-1,2t+1), Lz(N-2t-2,2t+2),...,Lz(M+1,M)$ provide the representations that differ from each other.
Therefore by Theorem $7$ following linear equations system are derived, if we regard $\Pi(X_{i}^{(2t+1)})$ as unknowns, $Lz(a,b)$ as coefficients\\
\begin{align*}
   & Lz(N-3,3)-c_3 Lz(N-1,1)= \sum p_{3j}\pi^{2j}Lz(a_{j},b_{j}) + \sum_{i=1}^{r_{3}}q_{3i}^{(3)}\Pi(X_{i}^{(3)})  \\
   & Lz(N-4,4)-c_4 Lz(N-1,1)= \sum p_{4j}\pi^{2j}Lz(a_{j},b_{j}) +  \sum_{i=1}^{r_{3}}q_{4i}^{(3)}\Pi(X_{i}^{(3)})  \\
    &Lz(N-5,5)-c_5 Lz(N-1,1)= \sum p_{5j}\pi^{2j}Lz(a_{j},b_{j}) +  \sum_{i=1}^{r_{3}}q_{5i}^{(3)}\Pi(X_{i}^{(3)}) +\sum_{i=1}^{r_{5}}q_{5i}^{(5)}\Pi(X_{i}^{(5)}) \\
   &...\\
   & Lz(M+1,M)-c_M Lz(N-1,1)= \sum p_{Mj}\pi^{2j}Lz(a_{j},b_{j}) + \sum_{i=1}^{r_{3}}q_{5i}^{(3)}\Pi(X_{i}^{(3)})+...+\sum_{i=1}^{r_{\tau}}q_{\tau i}^{(\tau)}\Pi(X_{i}^{(\tau)}) \\
\end{align*}
where $c,p,q\in \mathbb{Q}$, $a_{j}+b_{j}<N$, $\tau$ is the largest $2t+1$ such that $\mathcal{P}\mathcal{O}_{3}^{2t+1}(N)\neq \emptyset$.
There are $M-3+1$ equations. It obvious, the number of unknowns $|T(N)|=r_{3}+...+r_{\tau}\geq r_3$. Therefore if $M-3+1< r_{3}$, due to the unknowns are more than the number of equations, then $\Pi(X)=\prod_{(n,k)\in supp(X)}\zeta(n)^{k}$ cannot be solved by above equations system. It's well-known that $r_{3}=|\mathcal{P}_{\mathcal{O},3}^{3}(2M+1)|$ increases faster than $M$. Hence there exist $M_{1}$, such that if $M>M_{1}$, then $M_{1}-3+1< |\mathcal{P}_{\mathcal{O},3}^{3}(2M_{1}+1)|$. In fact, it's not hard to find out, if $M>9$, then $M-3+1<|T(2M+1)|$.\\
\begin{itemize}
  \item Case II. Suppose that $N=2M$ \\
\end{itemize}
Let\\
\begin{equation*}
  T(N)=\bigsqcup_{t=1}^{\infty}\mathcal{P}\mathcal{O}_{3}^{2t}(N)
\end{equation*}
Similar to the case of odd $N$, it reasonable to assume that if $Y\in \mathcal{P}_{2}(N)\backslash T(N)$, then
\begin{equation*}
  \sum Y=\sum p_{j} Lz(a_{j},b_{j})
\end{equation*}
with $a_{j}+b_{j}< N$. Now Suppose that \\
\begin{equation*}
\mathcal{P}\mathcal{O}_{3}^{2t}(N)=\{X_{1}^{(2t)},X_{2}^{(2t)},...,X_{r_{2t}}^{(2t)}\}
\end{equation*}
with $|\mathcal{P}\mathcal{O}_{3}^{2t}(N)|=r_{2t}$. According to theorem $7$, analogous linear equations system are derived
\begin{align*}
& Lz(N-2,2)-c_2 Lz(N-1,1)= \sum p_{2j}\pi^{2j}Lz(a_{j},b_{j}) + \sum_{i=1}^{r_{2}}q_{2i}^{(2)}\Pi(X_{i}^{(2)})  \\
& Lz(N-3,3)-c_3 Lz(N-1,1)= \sum p_{3j}\pi^{2j}Lz(a_{j},b_{j}) + \sum_{i=1}^{r_{2}}q_{3i}^{(2)}\Pi(X_{i}^{(2)})  \\
   & Lz(N-4,4)-c_4 Lz(N-1,1)= \sum p_{4j}\pi^{2j}Lz(a_{j},b_{j}) +  \sum_{i=1}^{r_{2}}q_{4i}^{(2)}\Pi(X_{i}^{(2)}) +\sum_{i=1}^{r_{4}}q_{4i}^{(4)}\Pi(X_{i}^{(4)}) \\
    &Lz(N-5,5)-c_5 Lz(N-1,1)= \sum p_{5j}\pi^{2j}Lz(a_{j},b_{j}) +  \sum_{i=1}^{r_{2}}q_{5i}^{(2)}\Pi(X_{i}^{(2)}) +\sum_{i=1}^{r_{4}}q_{5i}^{(4)}\Pi(X_{i}^{(4)}) \\
   &...\\
   & Lz(M,M)-c_M Lz(N-1,1)= \sum p_{Mj}\pi^{2j}Lz(a_{j},b_{j}) + \sum_{i=1}^{r_{2}}q_{5i}^{(2)}\Pi(X_{i}^{(2)})+...+\sum_{i=1}^{r_{\tau}}q_{\tau i}^{(\tau)}\Pi(X_{i}^{(\tau)}) \\
\end{align*}
where $c,p,q\in \mathbb{Q}$, $a_{j}+b_{j}<N$, $\tau$ is the largest $2t$ such that $\mathcal{P}\mathcal{O}_{3}^{2t}(N)\neq \emptyset$.
There are $M-2+1$ equations. It obvious, the number of unknowns $|T(N)|=r_{2}+...+r_{\tau}\geq r_2+r_4$ if $N$ is large enough. Therefore if $M-2+1< r_{2}+r_{4}$, due to the unknowns are more than the number of equations, then $\Pi(X)=\prod_{(n,k)\in supp(X)}\zeta(n)^{k}$ cannot be solved by above equations system. It's well-known that $r_{4}=|\mathcal{P}_{\mathcal{O},3}^{4}(2M)|$ increases faster than $M$. Hence there exist $M_{2}$, such that if $M>M_{2}$, then $M_{2}-2+1< |\mathcal{P}_{\mathcal{O},3}^{2}(2M_{2})|$. In fact, it's not hard to find out, if $M>10$, then $M-2+1<|T(2M)|$.\\

Finally, by above discussion we obtain $N_{0}=20$. If $N>N_{0}$, whenever $N$ is odd or even, there always exist $X\in \mathcal{P}_{2}(N)$ such that $X$ cannot be represented by finite $\mathbb{Q}(\pi)$-linear combination of $Lz(a,b)$ with $a+b\leq N$.
\begin{flushright}
  $\Box$
\end{flushright}
\end{pf}
\bibliographystyle{unsrt}  


Xiaowei Wang(\begin{CJK}{UTF8}{gbsn}王骁威\end{CJK})\\
Institut f\"{u}r Mathematik, Universit\"at  Potsdam, Potsdam OT Golm, Germany\\
 Email: \texttt{xiawang@gmx.de}
\end{document}